\documentclass[twoside,a4paper]{amsart}

\usepackage{mathabx}

\usepackage{url}      % Allows line-breaking of long URLs

\usepackage{pdfsync}
\usepackage{enumitem}
\usepackage{tikz}
\usetikzlibrary{calc}

\usepackage{xcolor}
\definecolor{webgreen}{rgb}{0,.5,0}
\definecolor{webbrown}{rgb}{.6,0,0}
\definecolor{RoyalBlue}{cmyk}{1, 0.50, 0, 0}
\usepackage[colorlinks=true, breaklinks=true, urlcolor=webbrown, linkcolor=RoyalBlue, citecolor=webgreen,backref=page]{hyperref}
\usepackage{bbold}

\usepackage{courier} % tt
\normalfont
\usepackage{epsfig, graphicx, subfigure}
\usepackage{amsmath, amssymb}

\newcommand{\T}     {\mathbb{T}}
\newcommand{\N}  {\mathbb{N}}
\newcommand{\D}     {\mathbb{D}}
\newcommand{\R}     {\mathbb{R}}

\newcommand{\Z}     {\mathbb{Z}}
\def\ge{\geqslant}
\def\le{\leqslant}

\newtheorem{theorem}{Theorem}[section]

\newtheorem{lemma}[theorem]{Lemma}

\numberwithin{equation}{section}

\begin{document}

\title[A sharp condition for the existence of wave operators \ldots  ]{
A sharp condition for the existence of wave operators   }

\begin{abstract}
We study the wave operators for a Jacobi matrix $J$ whose spectral measure $\rho$ satisfies the  Szeg\H{o} condition. For $\gamma_s$, the Verblunsky coefficients of the measure associated to $\rho$ by the Szeg\H o mapping, the condition $\lim_{n\to\infty}\log n\sum_{s=n}^{2n}|\gamma_s|^2=0$ guarantees (see \cite{DY}) the existence and completeness of wave operators.  In the current note, we prove that that condition is sharp. 
\end{abstract} \vspace{1cm}

%\thanks{
%The work of SD done in the first two sections
%  was supported by the grant NSF-DMS-1764245 and his research on the rest of the paper was supported by the Russian Science Foundation (project RScF-19-71-30004). The work of LM was supported by the grant RTG NSF-DMS-1147523.}

\author[S.A. Denisov]{Sergey A. Denisov}
\address{Department of Mathematics, University of Wisconsin-Madison, 480 Lincoln Dr., Madison, WI 53706, USA}
\email{\href{mailto:denissov@wisc.edu}{denissov@wisc.edu}}
\thanks{
This research was supported by the NSF grant DMS-2450716 and the Van Vleck Professorship Research Award. This work was also supported by  grants from  Simons Foundation International [SFI-MPS-SFM-00010669,SD] and [SFI-MPS-TSM-00012140,SD].  
}



\maketitle

\setcounter{tocdepth}{3}
\setcounter{tocdepth}{2}

\tableofcontents

%\section*{Contents}
%\begin{enumerate}
%\item Introduction
%\item Preliminary OPUC Results
%\item Existence and Completeness of Wave Operators
%\item Appendix
%\begin{enumerate}[label=\arabic{enumi}.\arabic*.]
%\item Auxiliary Proofs and an OPUC Estimate
%\item The Weak Nonlinear Carleson Condition
%\end{enumerate}
%\end{enumerate}

\section{Introduction and main result}

In \cite{DY}, we studied the existence of wave operators $\Omega_{\pm}$ for Jacobi matrices whose spectral measure belongs to the Szeg\H{o} class. In particular, we showed that the condition
\begin{equation}\label{sharp1}
 \lim_{n\to\infty}\log n\left(\sum_{s=n}^{2n}|\gamma_s|^2\right)=0
\end{equation}
on Verblunsky parameters $\{\gamma_n\}|_{n\in \Z^+}$ 
is sufficient for  the existence of $\Omega_{\pm}$.
In this note, we modify the construction in \cite{D2026} to prove that \eqref{sharp1} is sharp.

To set the stage,  we recall some notation from \cite{DY} making  changes when necessary.  Our introduction mostly follows \cite{DY}.  Consider semi-infinite Jacobi matrices 
\begin{align*}
J = \begin{pmatrix}
v_1& b_1&0&&&&\\
b_1& v_2&b_2&0&&&\\
0&b_2&v_3&b_3&0&&\\
&\ddots&\ddots&\ddots&\ddots&\ddots\\
\end{pmatrix}
\end{align*}
in the bounded and self-adjoint setting, where the Jacobi coefficients $\{ b_n\}$ and $\{v_n\}$ satisfy $b_n>0, v_n\in \R$ for all $n\in \N$, and $\{b_n\}, \{v_n\}\in \ell^\infty(\N)$.
Under these conditions, $J$ defines a bounded self-adjoint operator $J:\ell^2(\N)\to\ell^2(\N)$.  Like in \cite{DY}, our main interest  is the dynamics generated by $J$ through the discrete time-dependent Schr\"odinger equation:
\begin{align*}
-i\partial_t\psi=J\psi,\quad \psi(0)=f\in \ell^2(\N), 
\end{align*}
whose solution may be written in terms of the one-parameter unitary group as $\psi(t)=e^{itJ}f$.

%By Favard's theorem of classical analysis, there is a bijective correspondence between bounded semi-infinite Jacobi matrices $J$ and compactly and infinitely supported (in the sense of cardinality) probability measures on $\R$. In the forward direction of this correspondence, one maps $J$ to its canonical spectral measure, $\rho$. 

In this work, we consider Jacobi matrices whose corresponding measures $\rho$ are supported on $[-1,1]$, and are in the Szeg\H o class: given the decomposition of $\rho$ with respect to the equilibrium measure $\omega$ for the interval $[-1,1]$,
\begin{align*}
d\rho=\rho' d\omega+\rho_s,\,\,d\omega(x)=\pi^{-1}(1-x^2)^{-\frac 12}\chi_{[-1,1]}(x)dx,
\end{align*}
we assume that $\rho$ satisfies the Szeg\H o condition:
\begin{align}\label{eq:Szego}
\int\limits_\R \log \rho'd\omega>-\infty\,.
\end{align}
Denote by $J_0$ the free operator with Jacobi coefficients $b_n=\frac 12, v_n=0$ for all $n$. Like in \cite{DY}, our main objects of study are the wave operators for the pair $J$ and $J_0$. When they exist, these wave operators are defined as the strong $\ell^2(\N)$-limits
\begin{align}\label{eq:WO}
\lim_{T\to\pm \infty} e^{-iTJ}e^{iTJ_0}=:\Omega_{\pm}.
\end{align}
The question of whether $\Omega_{\pm}$ exist for Jacobi matrices in the Szeg\H{o} class was raised in \cite{DS}, pp. 33--34. In the current work, we  show that the Szeg\H{o} condition \eqref{eq:Szego} alone is not enough for the existence of wave operators. Additionally, we provide a sufficient condition for $\Omega_{\pm}$ to exist and show that it is sharp.\smallskip

%When $\Omega_{\pm}$ exist, the inclusion $\text{Ran}(\Omega_\pm)\subseteq \ell^2_{\rm ac}(\N)$ always holds. Completeness is defined as the equality $\text{Ran}(\Omega_\pm)=\ell^2_{\rm ac}(\N)$, so that $\Omega_\pm$ are unitary maps from $\ell^2(\N)$ onto $\ell^2_{\rm ac}(\N)$. This equality yields the comparison between the two evolutions alluded to above; indeed, once this is established, one has that every state in $\ell^2_{\rm ac}(\N)$ behaves asymptotically freely. Establishing these results under \eqref{eq:Szego}, along with a mild quantitative condition on the tail of the series for $\{ \gamma_n\}$, will be the main goal of our work here. For the reader more familiar with the OPUC literature, we note that this is a significantly more delicate task than existence and completeness of the CMV wave operators, where the free dynamics are given by translation, cf. \cite[Theorem 10.7.9]{SimonOPUC2}.\smallskip

%Despite the inverse direction of this correspondence passing through the orthogonal polynomials on the real line (OPRL), much of our work here makes use of orthogonal polynomials on the unit circle (OPUC). In particular, exploiting a well-known relationship between OPRL with respect to measures supported on the interval $[-1,1]$, and OPUC with respect to a related measure, we may work on the circle and apply harmonic analysis techniques. We refer the reader to \cite{SimonOPUC1,SimonOPUC2,simon2010szegHo,Szego} for further background on the brief description that follows.

Let $\sigma$ be an infinitely supported probability measure on the unit circle $\T$. Denote the orthonormal polynomials by $\phi_n(z)$ and the monic orthogonal polynomials by $\Phi_n(z)$.
Clearly, $\phi_n=\Phi_n/\|\Phi_n\|_{2,\sigma}$. We will use notation  $f^*(z)=z^n\overline{f(\bar z^{-1})}$ for a polynomial $f$ of degree $n$. The polynomials $\{\Phi_n\}|_{n=0}^\infty$ satisfy recurrence \cite[Theorem 1.5.2, p.~57]{SimonOPUC1}
\begin{equation}\label{opuc1}
\left\{\begin{array}{cc}
\Phi_{n+1}(z)=z\Phi_n(z)-\bar\gamma_n \Phi_n^*(z), &\Phi_0(z)=1,\\
\Phi_{n+1}^*(z)=\Phi_n^*(z)-z\gamma_n\Phi_n(z),&\Phi_0^*(z)=1\,,
\end{array}\right.
\end{equation}
where  $\gamma_n$ are the Verblunsky coefficients corresponding to the measure $\sigma$. The map $\sigma\mapsto \{\gamma_n\}|_{n=0}^\infty\in \D^\infty$ is a bijection, cf. \cite[p.~2]{SimonOPUC1}.

Measures $\rho$ supported on $[-1,1]$ may be related to even measures $\sigma$ on $\T$ through the Szeg\H o mapping $Sz$.  We will write $\sigma=Sz^{-1}(\rho)$ for the measure $\sigma$ on $\T$ that satisfies $d\sigma(z)=d\sigma(1/z)$ and
\begin{align}\label{opuc-oprl}
\int\limits_{\T}f(z)d\sigma(z)=\int\limits_{[-1,1]}f(e^{i\arccos(x)})d\rho(x)
\end{align}
for any $f\in C(\T)$ with $f(z)=f(1/z), z\in \T$ \cite[p.~31, eq.~(1.9.5)]{simon2010szegHo}. Throughout the text, $\sigma$ will refer to the measure related to $\rho$ in this way. For such measures, all $\gamma_n$ are real, and  they may be related to the Jacobi coefficients by means of the direct Geronimus relations \cite[p.~33]{simon2010szegHo}
\begin{align}\label{geronimus}
\begin{split}
b_{k+1}&=\frac 12 \left( (1-\gamma_{2k-1})(1-\gamma_{2k}^2)(1+\gamma_{2k+1})    \right)^{\frac 12},\\
v_{k+1}&=-\frac 12 (\gamma_{2k-2}(1+\gamma_{2k-1})-\gamma_{2k}(1-\gamma_{2k-1}))\,,
\end{split}
\end{align}
with the convention $\gamma_{-1}:= -1$. If we denote the polynomials orthonormal with respect to $\rho$ by $p_n(x)$, then \cite[p.~32, eq.~(1.9.12)]{simon2010szegHo}
\begin{align}\label{pnchin0}
p_n((z+1/z)/2)=(2(1-{\gamma}_{2n-1}))^{-1/2}z^{-n}(\phi_{2n}(z)+\phi_{2n}^*(z)), \, z\in \T\,.
\end{align}

Szeg\H o's theorem, cf. \cite[Theorem 1.8.6]{simon2010szegHo}, is a cornerstone of the theory of OPUC; for the decomposition of $\sigma$ with respect to the normalized Lebesgue measure $dm(z)$ on $\T$, $d\sigma(z)=\sigma'(z)dm(z)+\sigma_s(z)$, it states that
\begin{align}\label{nonp}
\prod_{n=0}^\infty (1-|\gamma_n|^2)=\exp\left(\int\limits_{\T}\log(\sigma'(z))dm(z) \right)\,,
\end{align}
which in turn implies
\begin{align*}
\{ \gamma_n\}\in \ell^2(\Z^+) \iff \int\limits_{\T} \log(\sigma'(z))dm(z)>-\infty.
\end{align*}
With this in hand, the identity \eqref{opuc-oprl} implies that the Szeg\H o condition \eqref{eq:Szego} on $J$ is equivalent to the corresponding Verblunsky parameters satisfying $\{\gamma_n\}\in \ell^2(\Z^+)$.

We denote by $\ell^2_{\rm ac}(\N)$ the absolutely continuous subspace for $J$, or set of vectors whose spectral measure is absolutely continuous, cf. \cite{LukicSpectralTheory,reed1972methods}:
\begin{align*}
\ell^2_{\rm ac}(\N):=\{\psi\in\ell^2(\N):\mu_{\psi} \text{ is a.c.} \},
\end{align*}
where $\mu_\psi$ is spectral measure of $\psi$ relative to $J$.
In \cite{DY}, working under the following additional assumption on the $\{\gamma_n\}$ related to $J$ as in the above, 
\begin{equation}\label{loga}
\lim_{n\to\infty} \log n\left(\sum_{s=n}^{2n}|\gamma_s|^2\right)=0,
\end{equation}
we proved 

\begin{theorem}[Denisov-Young, \cite{DY}]\label{thm:main} Given $\{\gamma_n\}\in \ell^2(\Z^+)$ and \eqref{loga}, the strong limits \eqref{eq:WO} exist.  Moreover, the wave operators are complete, i.e.,
\begin{align*}
{\rm Ran}(\Omega_\pm )=\ell^2_{\rm ac}(\N).
\end{align*}
\end{theorem}
In the current work, we show that condition \eqref{loga} is sharp. Specifically, the following statement holds
\begin{theorem}\label{mai}There is a sequence of  recurrence parameters $\{\gamma_n\}\in \ell^2(\Z^+)$ such that 
\begin{equation}\label{loga1}
 \limsup_{n\to\infty} \log n\left(\sum_{s=n}^{2n}|\gamma_s|^2\right)<\infty
\end{equation}
and the limits \eqref{eq:WO} do not exist.
\end{theorem}

%This work is organized as follows: the second section contains 

\textbf{Notation}
\begin{itemize}
%\item $J$ is a Jacobi matrix with Jacobi parameters $v_n$ on the diagonal and $b_n$ on the off-diagonal, and $\rho $ is its corresponding measure supported on $[-1,1]$.
%\item $J_0$ is the free Jacobi matrix with parameters $v_n=0$, $b_n=1/2$ corresponding to the measure $d\rho_0(x)=(1-x^2)d\omega(x)$ for $\omega$ is the equilibrium measure on $[-1,1]$.
\item $\Z^+=\{0,1,2,\dots\}$, $\N=\{1,2,3,\dots\}$, and $\R^+=[0,\infty)$.
\item $\rho$ is the canonical spectral measure for the Jacobi matrix $J$, and $\sigma$ is the even measure on $\T$ related to $\rho$ by the Szeg\H o mapping: $\rho=Sz(\sigma)$.
\item $\phi_n, \Phi_n$ for $n\in \Z^+$ are the orthonormal and monic orthogonal polynomials respectively with respect to $\sigma$, with $\{\gamma_n\}$ the Verblunsky parameters.
\item For a degree $n$ polynomial on the circle $f$, $f^*=z^n\overline{f}$. 
\item $p_n$, $n\in\Z^+$ are the orthonormal polynomials with respect to $\rho$.
\item $m$ is the  arclength measure on the torus $\T$ normalized to define a probability measure.
\item Given $\kappa\in \N$, we will denote the $\kappa$-th roots of unity by $z_{j,\kappa}:=\exp(2\pi ij/\kappa)$, for $j\in\{0,\dots,\kappa-1\}$. We refer to the argument as $x_{j,\kappa}:=2\pi j/\kappa$. 
\item Given a measure $\tau$, $\langle \cdot,\cdot\rangle_{\tau}$ stands for the inner product of two quantities in $L^2_{\tau}(\cdot)$.
\item $\|\cdot\|_{2,\tau}$ is the $L^2_\tau$ norm, and $\|\cdot\|_\infty$ is the $L^\infty$ norm, where the measure will be clear from context. 
%\item For a measure $\tau$ on $\T$ satisfying the Szeg\H o condition,
%\[
%\int\limits_\T \log(\tau'(z))dm(z)>-\infty,\,d\tau(z)=\tau'(z)dm(z)+\tau_s(z)
%\]
%the Szeg\H o function $D$ in the Hardy space $H^2(\D)$ is  \cite[p.~144, eq.~(2.4.2)]{SimonOPUC1}
%\begin{align}\label{eq:Szegofxn}
%D(z):=\exp\left(\frac{1}{2}\int\limits_\T\frac{w+z}{w-z}\log(\tau'(w))dm(w) \right)\,.
%\end{align}
%$D$ has nontangential boundary values $m$-almost everywhere so that we may write $D(z)$ when $|z|=1$ as well.
\item For a sequence $f=\{f_n\}\in \ell^2(\Z)$, $\widehat{f}\in L^2_{m}(\T)$ will denote the Fourier transform of $f$:
\[
\widehat{f}(z)=\sum_{n\in\Z}f_nz^n.
\]
For $g\in L^2_m(\T)$, its inverse is given by the sequence $\{g_n\}$ 
\[
\widecheck{g}(n)=g_n=\langle g,z^n\rangle_{m}=\int\limits_{\T}g(z)z^{-n}dm(z).
\]
For the function $f$ defined on $\T$, the symbol $\|f\|_{A_1(\T)}$ denotes the Wiener's norm
\[
\|f\|_{A_1(\T)}:=\sum_{n\in \Z} |f_n|\,.
\]
%\item For a function $f\in L^2(\R)$, the Fourier transform is denoted
%by
%\[
%(\mathcal{F}f)(\xi)=(2\pi)^{-1}\int\limits_\R f(x)e^{ -i x\xi}dx
%\]
%and the inverse Fourier transform is $\mathcal{F}^{-1}$.

\item Define  the $\mathcal{F}_{\rm non}: F\mapsto \frak{b}(z,F)$ where $F\in \ell^1(\Z)$ and $\frak{b}\in A_1(\T)$ denotes the second component of the SU(1,1) NLFT transform (see \cite{D2026}, p.2).  The symbol
$\mathcal{F}_{\rm non}^{-1}$ stands for the inverse map.

\item $\chi_E$ is the characteristic function for a set $E$, and $|E|$ denotes its Lebesgue measure. 
%\item For a sequence $\{a_j\}$ and a number $P$, $a_j^{(d)}=a_j\chi_{j\equiv d\mod P}$ for $0\leq d\leq P-1$. 
%\item ${\rm Proj}_Sv$  is the orthogonal projection of vector $v$ onto the subspace $S$ of a Hilbert space, and ${\rm Perp}_Sv=v-{\rm Proj}_Sv$ denotes its perpendicular component.
\item For quantities $A$ and $B$ and a parameter $r$, $A\le_{r}B$ means that there is a nonnegative constant $C_r$ with $A\leq C_r B$. $A\lesssim B$ means that $A\leq CB$ for a constant $C>0$ independent of the relevant variables, so that $A\lesssim 1$ is shorthand for $A$ being uniformly bounded in the relevant parameters. 
\item $A\sim B$ means $B\lesssim A\lesssim B$, and for a parameter $r$, $A\sim_r B$ means $B\le_{r} A\le_{r}B$.
%\item For a natural number $\ell$, the symbol $O_\ell(g)$ stands for a function $f$ that satisfies $|f|\le C_\ell g$ with some $\ell$-dependent constant $C_\ell$ and a nonnegative quantity $g$.
\item $\Omega_{\pm}$ are the strong limits \eqref{eq:WO}.
\end{itemize}\bigskip\bigskip

\section{Proof of the Theorem \ref{mai}}\label{section2}

\noindent {\it Proof of Theorem \ref{loga1}.} We first recall the ``daisy'' construction from \cite{D2026} and make the necessary modifications. Consider a large parameter $\kappa$ and the small positive parameters $\delta$ and $\eta$.
Later, we will choose a sequence of parameters $\kappa$ and $\delta$ as
\begin{equation}\label{paras}
\kappa_m:=2^{2^{2^m}}, \quad \delta_m:=\sqrt \eta(\log_2\kappa_m)^{-\tfrac 12}, \quad m\in \mathbb{N},\,\, m\ge m_0\gg 1,
\end{equation}
where $\eta$ is $m$-independent.

We recall the following definition from \cite{DY}. On the unit circle, we take a partition of unity $1=\sum_{j=0}^{\kappa-1} \omega_{j,\kappa}(z),\, \omega_{j,\kappa}(z):=\omega_0(z/z_{j,\kappa})$, $z\in \T$ where $z_{j,\kappa}=\exp(2\pi ij/\kappa), x_{j,\kappa}=2\pi j/\kappa$ and $\omega_{0}(e^{ix})=h(\kappa x)$ where $h$ is a suitable smooth nonnegative function with compact support in $[-2\pi,2\pi]$. Let \[
n_{j,\kappa}:=2\lceil \kappa^2\sin x_{j,\kappa}\rceil.\]

Using the perturbation theory for NLFT with small $\ell^1(\Z)$ norm, for every $j\in \{0,\ldots,\tfrac{\kappa}{10}\}$, we can find $b_{\kappa,\delta,j}: \T\to \R$ such that
\begin{eqnarray}\label{zhir}
\int b_{\kappa,\delta,j}^2dx\sim\frac{\delta^2}{\kappa}\,,\hspace{9cm}\\
\label{sok1}
|b_{\kappa,\delta,j}(e^{ix})|\lesssim \delta\Bigl((1+|\kappa (x-2\pi j/\kappa)|^2)^{-1}+(1+|\kappa (x+2\pi j/\kappa)|^2)^{-1}\Bigr), \quad |x|\le \pi,\\
|b'_{\kappa,\delta,j}(e^{ix})|\lesssim \delta\kappa\Bigl((1+|\kappa (x-2\pi j/\kappa)|^2)^{-1}+(1+|\kappa (x+2\pi j/\kappa)|^2)^{-1}\Bigr), \quad |x|\le \pi
\end{eqnarray}
and the sequence
\begin{equation}\label{sok2}
F_{\kappa,\delta,j}:=\mathcal{F}^{-1}_{\rm non}b_{\kappa,\delta,j}
\end{equation}
is real-valued, even, supported on $[-0.1\kappa,0.1\kappa]$ and satisfies
\begin{equation}\label{logs}
\|F_{\kappa,\delta,j}\|_{\ell^1(\Z)}\lesssim \delta\,.
\end{equation}
We sketch the construction of such a function later in the Appendix, see Lemma \ref{special}.
 Let $\frak{b}_{\kappa,\delta,j}:={b}_{\kappa,\delta,j}$ and $\frak{a}_{\kappa,\delta,j}^{(*)}$ be the associated outer function in $\D$ that satisfies
\[
|\frak{a}_{\kappa,\delta,j}^{(*)}(z)|^2=1+|\frak{b}_{\kappa,\delta,j}(z)|^2, \quad z\in \T\,.
\]
Clearly, $(\frak{a}_{\kappa,\delta,j},\frak{b}_{\kappa,\delta,j})$ is $SU(1,1)$ NLFT transform of the sequence $F_{\kappa,\delta,j}$.
Notice that the nonlinear Plancherel identity \eqref{nonp} gives
\begin{equation}\label{sok3}
\|F_{\kappa,\delta,j}\|_{\ell^2(\Z)}^2\sim \frac{\delta^2}{\kappa}\,.
\end{equation}
Let $\mu_{j,m}:=n_{j+\tfrac{\kappa_m}{8},\kappa_m}$. Notice that $\mu_{j+1,m}-\mu_{j,m}\ge \kappa_m/4$ as long as $j\in \{0,\ldots,\tfrac{\kappa_m}{16}\}.$
The single $\kappa$ and $\delta$ ``daisy'' pattern in \cite{D2026} is now modified as follows 
\begin{equation}\label{stop1}
Q_m(n)=\sum_{1\le j\le \tfrac{\kappa_m}{16}} F_{\kappa_m,\delta_m,j}(n-\mu_{j,m})
\end{equation}
and we choose $F$, the NLFT sequence, as an infinite sum of well-separated ``daisies''
\begin{equation}\label{stop2}
F(n)=\sum_{m\geq m_0} Q_{m}(n)
\end{equation}
with parameters $\kappa_m$ and $\delta_m$ defined in \eqref{paras}.
By construction, $F$ is supported on $\Z^+$.
From \eqref{sok3}, we get
\begin{equation}\label{pes1}
\|Q_{m}\|_{\ell^2}\sim \delta_m^2, \quad
\|F\|^2_{\ell^2}\sim \sum_{m\ge m_0}\delta_m^2=\eta \sum_{m\ge m_0}2^{-2^m}\sim \eta 2^{-2^m}\,.
\end{equation}
We now define Verblunsky parameters of matrix $J$  by  letting
\[
\gamma_n:=F(n), \quad n\ge 0\,.
\]
As before,  denote the orthogonality measure by $\sigma$ and corresponding orthogonal polynomials by $\Phi_n(z)$ and $\phi_n(z)$.
First, we have $\gamma_n\in \R$ and $\|\{\gamma_n\}\|_{\ell^2(\Z^+)}\ll 1$. Next, we check that
\[
\log_2\mu_{0,m}\sum_{j=\mu_{0,m}}^{C\mu_{0,m}} |\gamma_j|^2 \stackrel{\eqref{pes1}}{\sim} (\log_2\kappa_m)\cdot  \delta_{m}^2\sim \eta
\]
with any fixed constant $C$.
Due to sparseness of the sequence $\{\kappa_m\}$, the bound
 \eqref{loga1} holds with a positive upper limit. 

Now, we turn to the analysis of the NLFT of truncations, i.e., $\frak{a}(F\cdot \chi_{n\le N})$ and $\frak{b}(F\cdot \chi_{n\le N})$ for different $N$. Since the shifted $F_{\kappa_m,\delta_m,j}$ in the sums in \eqref{stop1} and \eqref{stop2} all have disjoint compact supports of size at most $0.2\kappa_m$, we can apply Lemma 2.3 from \cite{D2026} to study $\frak{a}(F\cdot \chi_{n\le N})$ and $\frak{b}(F\cdot \chi_{n\le N})$.  Let $\nu_{j,m}=\mu_{j,m}+\tfrac{\kappa_m}{10}$. Within each daisy, the recursion gives for $j\in \{1,\ldots,\tfrac{\kappa_m}{16}\}$:
\[
\begin{array}{c}
|\frak{a}(z,Q_m\cdot \chi_{n\le\nu_{j,m}})|\le |\frak{a}_{\kappa_m,\delta_m,j}(z)|\cdot|\frak{a}(z,Q_{m}\cdot \chi_{n\le \nu_{j-1,m}})|+|\frak{b}_{\kappa_m,\delta_m,j}(z)|\cdot|\frak{b}(z,Q_{m}\cdot \chi_{n\le \nu_{j-1,m}})|,\\
|\frak{b}(z,Q_m\cdot \chi_{n\le\nu_{j,m}})|\le |\frak{b}_{\kappa_m,\delta_m,j}(z)|\cdot|\frak{a}(z,Q_{m}\cdot \chi_{n\le \nu_{j-1,m}})|+|\frak{a}_{\kappa_m,\delta_m,j}(z)|\cdot|\frak{b}(z,Q_{m}\cdot \chi_{n\le \nu_{j-1,m}})|.
\end{array}
\]
Restrict $z=e^{ix}, x\in [0,\pi]$ and use \eqref{sok1}  along with an adaptation of Lemma \ref{adjust} from the Appendix to get
\[
|\frak{b}(e^{ix},Q_m\cdot \chi_{n\le  \nu_{j,m} })|\lesssim \delta_m \left\{ 
\begin{array}{cc}
1, & 0\le x\le 2\pi j/\kappa_m,\\
(1+|\kappa_m(2\pi j/\kappa_m-x)|)^{-1}, & 2\pi j/\kappa_m<x<\pi
\end{array}
\right.
\]
and
\begin{equation}\label{puk1}
\frak{a}^{(*)}(e^{ix},Q_{m}\cdot \chi_{n\le \nu_{j,m}})=\left(\prod_{s=1}^j\frak{a}^{(*)}_{\kappa_m,\delta_m,s}(z)\right)(1+O(\delta_m))
\end{equation}
uniformly over $z=e^{ix}\in \T$. Then, 
\[
\left|\text{arg} \left(\prod_{s=1}^j\frak{a}^{(*)}_{\kappa_m,\delta_m,s}(e^{ix})\right)\right| \sim \sum_{p=-j}^{j} \frac{\delta_m^2}{1+\kappa_m|x-\tfrac{2\pi p}{\kappa_m}|}, \quad 2\pi j/\kappa_m<x<\pi\,.
\]
Hence, using \eqref{puk1}, we get for $j>\kappa_m/100$
\[
\left|\text{arg}\,\frak{a}^{(*)}(e^{ix},Q_{m}\cdot \chi_{n\le \nu_{j,m}})\right| \sim  {\delta_m^2\log\kappa_m}\sim \eta, \quad |x-2\pi j/\kappa_m|<C\kappa_m^{-1}
\]
with any fixed $C$. Similarly, the recursions that define the NLFT give
\[
\left|\text{arg}\,\frak{a}^{(*)}(e^{ix},F\cdot \chi_{n\le \nu_{j,m}})-\text{arg}\,\frak{a}^{(*)}(e^{ix},F\cdot \chi_{n\le \nu_{0,m}})\right| \sim  {\delta_m^2\log\kappa_m}\sim \eta, \quad |x-2\pi j/\kappa_m|<C\kappa_m^{-1}\,.
\]
Now, we recall that (see (3.2) in \cite{D2026})
\[
\phi^*_N(z)=\frak{a}^{(*)}(z,F\cdot\chi_{n\le N})+\frak{b}(z,F\cdot\chi_{n\le N})=\frak{a}^{(*)}(z,F\cdot\chi_{n\le N})(1+\frak{r}_N)
\]
where
\[
\frak{r}_N:=\frak{b}(z,F\cdot\chi_{n\le N})/\frak{a}^{(*)}(z,F\cdot\chi_{n\le N})\,.
\]
In \cite{D2026}, we saw that $\frak{r}_N$ converges to $\frak{r}$ uniformly on $\T$ as $N\to\infty$ and that $\|\frak{r}\|_{L^\infty(\T)}=o(1), m_0\to\infty$ and 
$\sup_{N\in \N}\||\phi^*_N(e^{ix})|-1\|_{L^\infty(\T)}=o(1),$ $m_0\to\infty$. We also showed in \cite{D2026} that the orthogonality measure $\sigma$ is given by a positive and continuous weight. Then, for sufficiently large $m_0$, we get
\begin{equation}\label{dva-sd3}
|\text{arg} \,\phi_{\mu_{j,m}}^*(e^{ix})-\text{arg} \,\phi_{\mu_{0,m}}^*(e^{ix})|\sim \eta
\end{equation}
for $|x-2\pi j/\kappa_m|<C\kappa_m^{-1}$. Therefore, we get 
\begin{equation}\label{dva-sd4}
|\,\phi_{\nu_{j,m}}^*- \,\phi_{\nu_{0,m}}^*|\sim \eta
\end{equation}
for $|x-2\pi j/\kappa_m|<C\kappa_m^{-1}$ and $\tfrac{\kappa_m}{16}>j>\tfrac{\kappa_m}{100}$. From \eqref{logs}, one has
\begin{equation}\label{dva-sd5}
|\,\phi_{\mu_{j,m}}^*- \,\phi_{\mu_{0,m}}^*|\sim \eta
\end{equation}
when $|x-2\pi j/\kappa_m|<C\kappa_m^{-1}$ and $\tfrac{\kappa_m}{16}>j>\tfrac{\kappa_m}{100}$.

\bigskip

In \cite{DY}, we showed that existence of wave operators $\Omega_{\pm}$ implies (see formula (3.11) in \cite{DY}) that
\[
\sum_{1\le j\le \kappa_m/16}  (\phi^*_{\mu_{j,m}}-\phi^*_{\mu_{0,m}})\omega_{j+\tfrac{\kappa_m}{8},\kappa_m}\to 0
\]
 in $L^2_\sigma(I)$ for any interval $I\subset (0,\pi/2)$. That, however, contradicts \eqref{dva-sd5}. \qed

\bigskip

\section{Appendix}

Suppose $b(x)$ is a nonnegative function defined on $\R$ that satisfies 
\begin{equation}
b(x)\lesssim \frac{\delta}{1+|x|^{\gamma}}, \quad \gamma>1, \quad \delta<1.
\end{equation}
and let $a(x)$ satisfy
\begin{equation}
1\le a(x), \quad a^2(x)=b^2(x)+1\,.
\end{equation}
\begin{lemma}\label{adjust}
Suppose $A_n(x)$ and $B_n(x)$ are nonnegative function defined on $\R$ that satisfy recursions for $n\in \Z^+$:
\begin{equation}\left\{
\begin{matrix}
A_{n+1}(x)&\le& a(x-n)A_n(x)+b(x-n)B_n(x), \\
B_{n+1}(x)&\le& b(x-n)A_n(x)+a(x-n)B_n(x),
\end{matrix}\right.
\end{equation}
and $A_0(x)=a(x),\, B_0(x)=b(x)$. Then, 
\begin{equation}\label{telescop}
B_n(x)\lesssim_\gamma \delta\sum_{j=0}^n b(x-j)\lesssim_\gamma \delta\left\{\begin{array}{cc}
(1+|x-n|)^{-\gamma+1},& x>n,\\
1,& 0\le x\le n\\
\end{array}\right.\,.
\end{equation}
\end{lemma}
\begin{proof}
Let
\[
\Omega_n(x)=\prod_{j=0}^n a(x-j),\, A_n=\Omega_n \widehat{A}_n,\, B_n=\Omega_n\widehat{B}_n\,.
\]
So, 
\begin{equation}\label{skok1}
\left\{
\begin{matrix}
\widehat A_{n+1}(x)&\le & \widehat A_n(x)+\frac{b(x-n)}{a(x-n)}\widehat B_n(x), \quad \widehat A_0(x)=1\,,\\
\widehat B_{n+1}(x)&\le & \frac{b(x-n)}{a(x-n)}\widehat A_n(x)+\widehat B_n(x),\quad  \widehat B_0(x)=b(x)/a(x)\,.
\end{matrix}\right.
\end{equation}
Since
\[
\left\|\begin{matrix}
1& \frac{b(x-n)}{a(x-n)} \\
\frac{b(x-n)}{a(x-n)}&1
\end{matrix}\right\|=1+\frac{b(x-n)}{a(x-n)}, \quad \sum_n \frac{b(x-n)}{a(x-n)}\le \sum_n {b(x-n)}\lesssim_\gamma \delta
\]
we immediately get $\widehat A_n(x)\lesssim_\gamma 1$ uniformly in $x$ and $n$. Hence, the second bound in \eqref{skok1} gives us 
\[
\widehat{B}_n(x)\lesssim_\gamma \sum_{j=0}^n b(x-j)\,.
\]
Next, we notice that $1\le |\Omega_n|\lesssim_\gamma 1$ uniformly in $n$ and $x$. We combine it with the previous estimate to  get \eqref{telescop}.
\end{proof}

We now state and sketch the proof of the construction of the bump functions required in the proof of the main result.
Take $\kappa\in \N, \delta\in (0,0.1)$ and $j\in \{0,\ldots,\lceil\kappa/10\rceil\}$. Let $K_\kappa(z)$ be a scaled Fejer kernel of degree $\kappa$:
\[
K_{\kappa,\delta}(e^{ix})=\frac{\delta}{\kappa^2}\left(\frac{\sin(\kappa x/2)}{\sin(x/2)}\right)^2
\]
so the formula $\widecheck{K}_{\kappa,\delta}(n)=\delta\kappa^{-1}(1-|n|/\kappa)\cdot \chi_{|n|\le \kappa}$ holds for  its inverse Fourier transform. In inverse scattering theory, we can address the function $\mathcal{F}_{non}(\widecheck{K}_{\kappa,\delta}(n))$ as {\it the nonlinear scaled Fejer kernel} and the perturbation theory can be used to show that this function is close to $K_{\kappa,\delta}(e^{ix})$ for small $\delta$.
\begin{lemma}\label{special}
Define $b_{\kappa,\delta,j}(e^{ix})=\mathcal{F}_{non}(\widecheck{K}_{\kappa,\delta}(n) \cdot 2\cos(2\pi j n/\kappa))$. Then, $b_{\kappa,\delta,j}$ is a real-valued trigonometric polynomials of degree at most $\kappa$ that satisfies
\begin{eqnarray}\label{zhir1}
\int b_{\kappa,\delta,j}^2dx\sim\frac{\delta^2}{\kappa}\,,\hspace{9cm}\\
\label{sok11}
|b_{\kappa,\delta,j}(e^{ix})|\lesssim \delta\Bigl((1+|\kappa (x-2\pi j/\kappa)|^2)^{-1}+(1+|\kappa (x+2\pi j/\kappa)|^2)^{-1}\Bigr), \quad |x|\le \pi,\\
|b'_{\kappa,\delta,j}(e^{ix})|\lesssim \delta\kappa\Bigl((1+|\kappa (x-2\pi j/\kappa)|^2)^{-1}+(1+|\kappa (x+2\pi j/\kappa)|^2)^{-1}\Bigr), \quad |x|\le \pi\,.\label{sok12}
\end{eqnarray}

\end{lemma}
\noindent{\it Sketch of the proof.} The proof is a routine  argument in perturbation theory for $\ell^1(\Z)$ space of coefficients.
The iteration of the NLFT recursion gives the following representation for $b_{\kappa,\delta,j}$:
\begin{equation}\label{kot}
b_{\kappa,\delta,j}(e^{ix})=C_{\kappa,\delta,j} \sum_{r \text{\,is \.odd}}\sum_{n_1<\ldots<n_r}q(n_1)\cdot\ldots\cdot q(n_r)e^{ix(n_1-n_2+\ldots+n_r)}\,,
\end{equation}
where we set $q(n)=2\widehat{K}_\kappa(n)\cos(2\pi j n/\kappa)$ for shorthand. We have $C_{\kappa,\delta,j}\sim 1$ uniformly in the parameters.
The formula \eqref{kot} shows that $\deg b_{\kappa,\delta,j}\le \kappa$ and that it is real-valued. We also have
\begin{eqnarray*}
\left\|\sum_{n_1<\ldots<n_r}q(n_1)\cdot\ldots\cdot q(n_r)e^{ix(n_1-n_2+\ldots+n_r)}\right\|_{A_1(\T)}\lesssim \|q\|_{\ell^1(\Z)}^r\le (4\delta)^r\,,\\ \left\|\sum_{n_1<\ldots<n_r}q(n_1)\cdot\ldots\cdot q(n_r)e^{ix(n_1-n_2+\ldots+n_r)}\right\|_{L^\infty(\T)}\lesssim \|q\|_{\ell^1(\Z)}^r/r!\le (4\delta)^r/r!\,.
\end{eqnarray*}
The bound \eqref{zhir1} follows from the nonlinear Plancherel identity \eqref{nonp}. The estimate \eqref{sok12} is the corollary of \eqref{sok11} and the convolution representation for the trigonometric polynomial of degree at most $\kappa$. Hence, we only need to focus on proving \eqref{sok11}. That estimate, however, can be obtained by the standard summation by parts (twice) in each of the $r$-term and by using the flatness of $\widehat K_{\kappa,\delta}$ in variable $n$. For example, for $r=1$, we use an identity
\[
(1-z)^2\sum_{n_1}\widecheck K_{\kappa,\delta}(n_1)z^{n_1}=\sum_{n_1}(\widecheck K_{\kappa,\delta}(n_1)+\widecheck K_{\kappa,\delta}(n_1-2)-2\widecheck K_{\kappa,\delta}(n_1-1))z^{n_1}
\]
which gives
\[
|(1-z)^2\sum_{n_1}\widecheck K_{\kappa,\delta}(n_1)z^{n_1}|=|\sum_{n_1}(\widecheck K_{\kappa,\delta}(n_1)+\widecheck K_{\kappa,\delta}(n_1-2)-2\widecheck K_{\kappa,\delta}(n_1-1))z^{n_1}|\lesssim \delta\kappa^{-2}\,.
\]
Coupled with an estimate
\[
|\sum_{n_1}\widecheck K_{\kappa,\delta}(n_1)z^{n_1}|\lesssim \delta\,,
\]
we have
\[
\left|\sum_{n_1}\widecheck K_{\kappa,\delta}(n_1)e^{in_1x}\right|\lesssim \delta/(1+(\kappa x)^2)\,.
\]
Substituting  $q$ instead of $\widecheck K_{\kappa,\delta}$, we only create two shifts in $x$-variable, which results in the bound by the right hand side of \eqref{sok11}. The terms corresponding to $r>1$ can be handled similarly.
\qed

\bigskip
\noindent {\bf AI usage.} No AI tools were used in preparing that paper except for the verification of routine calculations in Lemma \ref{special} for which the ChatGPT-6 Astra was utilized.

\end{document}